\input amstex
\documentstyle{amsppt}
%\NoRunningHeads\NoPageNumbers\nologo
\footline{}

\catcode`\@=11
\def\raggedcenter@{\leftskip\z@ plus.3\hsize \rightskip\leftskip
 \parfillskip\z@ \parindent\z@ \spaceskip.3333em \xspaceskip.5em
 \pretolerance9999\tolerance9999 \exhyphenpenalty\@M
 \hyphenpenalty\@M \let\\\linebreak}
\catcode`\@=\active

\mag=1200
\pagewidth{14truecm}
\pageheight{21.5truecm}
\hoffset23pt
\voffset36pt
\binoppenalty=10000
\relpenalty=10000
\tolerance=500
\mathsurround=1pt

\define\m1{^{-1}}
\define\gp#1{\langle#1\rangle}
\define\bk{1}
\define\bu{2}

\hrule height0pt
\vskip 1truein

\topmatter

\title
Symmetric units in modular group algebras
\endtitle

\author
Victor Bovdi\\
{\eightpoint Bessenyei Teachers College,
4401 Ny\'\i regyh\'aza, Hungary}\\
\phantom{anything}\\
 L.\,G. Kov\'acs\\
{\eightpoint Australian National University,
Canberra, ACT 0200, Australia}\\
\phantom{anything}\\
and\\
\phantom{anything}\\
S.\,K. Sehgal\\
{\eightpoint University of Alberta, Edmonton, Alberta T6G 2H1, Canada}
\endauthor

\leftheadtext\nofrills{Victor Bovdi, L.\,G. Kov\'acs and S.\,K. Sehgal}
\rightheadtext\nofrills{Symmetric units in  group algebras}

\thanks
V. Bovdi is indebted to the University of Alberta for warm hospitality and
generous support
during a period when part of this work was done.  His research was also
supported by the
Hungarian  National Foundation for Scientific Research grant No.\,T 014279.
\endthanks

\abstract\baselineskip12pt
Let $p$ be a prime, $G$ a locally finite $p$-group, $K$ a commutative ring
of characteristic
$p$.  The anti-automorphism $g\mapsto g\m1$ of $G$ extends linearly to an
anti-automorphism $a\mapsto a^*$ of $KG$. An element $a$ of $KG$ is called
symmetric if
$a^*=a$. In this paper we answer the question: for which $G$ and $K$ do the
symmetric
units of $KG$ form a multiplicative group.
\endabstract


\endtopmatter

\addto\tenpoint{\normalbaselineskip16pt\normalbaselines}

\document

Let $G$ be a group, $K$ a commutative ring (with $1$), and $U(KG)$ the
group of units in the
group algebra $KG$. The anti-automorphism $g\mapsto g\m1$ of $G$ extends
linearly to an
anti-automorphism $a\mapsto a^*$ of $KG$; this extension leaves $U(KG)$ setwise
invariant. An element $a$ of $KG$ is called {\it symmetric\/} if $a^*=a$.

It is an open problem to find the noncommutative $KG$ in which the
symmetric units form a
multiplicative group. Here we solve this under the assumption that $K$ has prime
characteristic, $p$ say, and $G$ is a locally finite $p$-group. Our result
is the following.

\proclaim
 {Theorem} Let $p$ be a prime, $K$ a commutative ring of characteristic $p$, and
$G$ a nonabelian locally finite $p$-group. The symmetric units of $KG$ form
a multiplicative
group if and only if $p=2$ and $G$ is the direct product of an elementary
abelian group and a
group $H$ for which one of the following holds:
\itemitem{\bf(i)} $H$ has an abelian subgroup $A$ of index $2$ and an
element $b$ of order
$4$ such that conjugation by $b$ inverts each element of $A$;
\itemitem{\bf(ii)} $H$ is the direct product of a quaternion group of order
$8$ and a cyclic
group of order $4$, or the direct product of two quaternion groups of
order~$8$;
\itemitem{\bf(iii)} $H$ is the central product of the group
$\gp{\,x,y\mid x^4=y^4=1,\ x^2=[y,x]\,}$  with a quaternion group of order $8$,
the nontrivial element common to the two central factors being $x^2y^2$;
\itemitem{\bf(iv)} $H$ is isomorphic to one of the groups $H_{32}$ and
$H_{245}$ defined
below.
\endproclaim

The relevant definitions are:
$$\align
 H_{32}=\bigl\langle\,x,y,u\bigm|\ &x^4=y^4=1,\\
&x^2=[y,x],\\
&y^2=u^2=[u,x],\\
&x^2y^2=[u,y]\,\bigr\rangle,\\
\noalign{\vskip2pt}\allowbreak
 H_{245}=\bigl\langle\,x,y,u,v\bigm|\ &x^4=y^4=[v,u]=1,\\
&x^2=v^2=[y,x]=[v,y],\\
&y^2=u^2=[u,x],\\
&x^2y^2=[u,y]=[v,x]\,\bigr\rangle.\\
\endalign
$$

Note that in case (i) all elements of $H$ outside $A$ have order $4$ and so
any one of them
can serve as $b$. The list of groups in this theorem is part of the list in Theorem 1.2 of Bovdi
and Kov\'acs \cite{\bk}, and the proof relies heavily on Lemma 1.4 of that
paper.

\medskip

The proof of the Theorem will occupy the rest of this note.

Set $S=\{\,t\mid t\in G,\,t^2=1\,\}\cup\{\,g+g\m1\mid g\in
G,\,g^2\neq1\,\}$, and
note that the symmetric elements of $KG$ are precisely the $K$-linear
combinations of the
elements of $S$.  Like the fixed points of any anti-automorphism of any
group, the symmetric
units form a subgroup in $U(KG)$ if and only if they commute with each
other.  It is well known
that once $K$ is of characteristic $p$ and $G$ is a locally finite
$p$-group, the augmentation
ideal of $KG$ is locally nilpotent, and so every element congruent to $1$
modulo this ideal is a
unit. In particular, $-1+g+g\m1$ is a symmetric unit in $KG$ whenever $g\in
G$. Of course,
$t$ is a symmetric unit whenever $t\in G$ and $t^2=1$. This proves that
{\it the symmetric
units form a multiplicative group if and only if every pair of elements of
$S$ commutes.}

In particular, for a given $p$ this issue is independent of the choice of
$K$. It will be
convenient to call a locally finite $p$-group {\sl good}  if every pair of
elements in $S$
commutes (say, in $KG$ with $K=\Bbb Z/p\Bbb Z)$. Note that all abelian
groups are
good, all subgroups of good groups are good, and that a locally finite
$p$-group is good if all its
$2$-generator subgroups are good.

In a good group, any two involutions (that is, elements of order $2$)
commute. If $g$
and $t$ are as in the definition of $S$, then the only way $t$ and $g+g\m1$
can commute is if
$g^t\ (=t\m1gt)$ is either $g$ or $g\m1$. In the second case the subgroup
$\gp{g,h}$
generated by $g$ and $t$ is a nonabelian dihedral group and therefore
contains noncommuting
involutions. This proves that {\it in a good group every involution is
central.}

Next we prove that if $g$, $h$ are noncommuting elements in a good group
$G$, then there
exist $x$, $y$ in $G$ such that $\gp{g,h}=\gp{x,y}$ and $x^y=x\m1$. To this
end, note that
any two of $g$, $h$, $gh$ generate the nonabelian group $\gp{g,h}$, so by
the previous
paragraph none of them can have square $1$. On the other hand, $g+g\m1$ and
$h+h\m1$
commute: thus
$$
gh+gh\m1+g\m1h+g\m1h\m1=hg+hg\m1+h\m1g+h\m1g\m1.
$$
 If $gh$ occurs more than once on the left hand side, we must have
$gh=g\m1h\m1$\!, so
$x=gh$, $y=h$ will do. Otherwise $gh$ must equal one of the summands on the
right hand side.
That summand cannot be $hg$, for $g$ and $h$ do not commute; nor can it be
$h\m1g\m1$,
for $(gh)^2\neq1$. Thus either $gh=hg\m1$, in which case $x=g$, $y=h$ will
work, or
$gh=h\m1g$, and then we can take $x=h$, $y=g$.

This already shows that the prime $p$ involved in a nonabelian good group
can only be $2$.

The point established in the second last paragraph can be taken further: in
those circumstances,
$x$ and $y$ can be chosen so that the order of $y$ is~$4$. To see this,
note first that
$(xy)\m1=xy\m1\neq xy$,  hence $xy+xy\m1$ and $y+y\m1$ must commute. Given that
the characteristic is $2$, this leads to $x(y^2+y^{-2})=x\m1(y^2+y^{-2})$.
If the
cosets $x\gp{y}$ and $x\m1\gp{y}$ are different, this forces
$y^2+y^{-2}=0$, that is,
$y^4=1$. If $x\gp{y}=x\m1\gp{y}$, it is easy to deduce that this is the
only nontrivial coset
of $\gp{y}$ in $\gp{x,y}$. Groups of $2$\hbox{-}power order with a cyclic
subgroup of index
$2$ are well known (see for example Section 109 in Burnside's book
\cite{\bu}); the
nonabelian groups of this kind with all involutions central are precisely
the generalized
quaternion groups. (We count the quaternion group of order $8$ among the
generalized
quaternion groups.) Of course, each generalized quaternion group can be
generated by a pair of
elements $x$, $y$ such that $x^y=x\m1$ and $y^4=1$.

We sum up these conclusions in the following.

\proclaim
 {Lemma 1}
 If $G$ is a nonabelian good group, then $p=2$ and each nonabelian $2$-generator
subgroup of $G$ is either a generalized quaternion group or a semidirect
product
$$
C_{2^m}\rtimes C_4=\bigl\langle\,x,y\bigm|x^{2^m}=y^4=1,\ x^y=x\m1\,\bigr\rangle
$$
with $m\geq2$.\qed
\endproclaim

If $m>2$ then there are other semidirect products that might be called
$C_{2^m}\rtimes C_4$, but we shall always mean this one.

\proclaim
 {Lemma 2}
 If $G$ is a nonabelian good group and the exponent of $G$ is not $4$, then
$G$ has an
abelian subgroup $A$ of index $2$ and an element $b$ of order $4$ such that
conjugation by
$b$ inverts each element of $A$.
\endproclaim

\demo{Proof}
  Set $A=\bigl\langle\,a\in G\bigm|a^4\neq1\,\bigr\rangle$. Suppose first
that $A$ is
nonabelian. Then there are noncommuting elements $g$, $h$ in $A$ whose
orders are greater
than~$4$. In a generalized quaternion group, all elements of order greater
than~$4$ lie in one
cyclic subgroup, so $\gp{g,h}$ cannot be a generalized quaternion group. In a
$C_{2^m}\rtimes C_4$, all elements outside $\gp{x,y^2}$ have order~$4$, so
$\gp{g,h}$ cannot be a $C_{2^m}\rtimes C_4$ either. This contradiction to
Lemma~1
proves that $A$ must be abelian.

Let $b$ be any element of $G$ outside $A$: by the definition of $A$, then
$b^4=1$. If
$a^4\neq1$, then $a$ and $b$ cannot commute (else we would have
$(ab)^4\neq1$ and then
$a,\, ab\in A$,\ $b\in A$ would follow, contrary to the choice of $b$). In
a generalized
quaternion group or in a $C_{2^m}\rtimes C_4$, an element of order greater
than~$4$ can
only be conjugate to itself or to its inverse, so Lemma~1 implies that
$a^b=a\m1$. It
follows that $b$ inverts every element of $A$. That includes $b^2$, so
$b^4=1$, and as all
involutions are central we cannot have $b^2=1$.\qed
\enddemo

It follows from Lemma 1 that if $G$ is a nonabelian good group of exponent
$4$ and if
$g,\,h\in G$, then $\gp{g^2}$ is central and $\gp{g,h}/\gp{g^2}$ is abelian
or dihedral.
Under somewhat weaker hypotheses, Lemma 1.4 of \cite{\bk} asserts that $G$
is a direct
product of an elementary abelian group and a group $H$ such that either $H$
satisfies one of
the conditions (i)--(iv) of our Theorem, or $H$ is an extraspecial
$2$-group, or $H$ is the
central product of an extraspecial $2$\hbox{-}group with a cyclic group of
order $4$. All
central products of this kind and all extraspecial $2$-groups except the
quaternion group
contain noncentral involutions, while the quaternion group satisfies
condition~(i). This
completes the proof of the `only if' part of our Theorem.

\medskip

The proof of the `if' part is much easier. The definition of good group
directly yields that the
direct product of an elementary abelian $2$-group and a good $2$-group is
always good: thus
it suffices to check that if $p=2$ and one of the conditions (i)--(iv)
holds then $H$ is
good.

Consider case (i) first. We have already remarked that in this case all
involutions of $H$ lie in
$A$, so they are all central. If $g\in A$, then $g+g\m1$ commutes with
every element of
$H$ and is therefore central in $KH$. If both $g$ and $h$ are elements of
$H$ outside $A$,
one can play the role of $b$ and the other can be written as $ab$ with
$a\in A$, and
$(b+b\m1)^2=0$ implies that
$(b+b\m1)((ab)+(ab)\m1)=0=((ab)+(ab)\m1)(b+b\m1)$. Thus
in this case $H$ is good.

In the other three cases $H$ has exponent $4$ and we know (from O'Brien's
Lemma 4.1 in
\cite{\bk}, or by direct calculation) that all involutions in $H$ are
central. We can also
see that the Frattini subgroup $\gp{\,h^2\mid h\in H\,}$ of $H$ has
order~$4$. Thus
if $g,\,h\in H$, then $\gp{g,h}$ has order at most~$16$. The groups of
order dividing~$16$
are well known (see for example Section 118 in~\cite{\bu}); there are only two
$2$-generator nonabelian groups of exponent~$4$ among them in which all
involutions are
central, and both of those satisfy condition~(i).  Thus by the previous
paragraph $\gp{g,h}$ is
good, and so $g+g\m1$ commutes with $h+h\m1$. We conclude that $H$ is good
in each of
these cases.

This completes the proof of the Theorem.

\enddocument